\documentclass[a4paper,10pt,reqno]{amsart}
\usepackage[a4paper,tmargin=40mm,bmargin=35mm,hmargin=30mm]{geometry}
\usepackage[T1]{fontenc}
\usepackage{lmodern}
\usepackage{amsmath}
\usepackage{amssymb}
\usepackage{amsthm}
\usepackage{stmaryrd}
\usepackage{tikz}
\usepackage[all]{xy}
\usepackage{booktabs}
\usepackage{multirow}
\usepackage{multicol}
\usepackage{here}
\usepackage{aliascnt}
\usepackage{hyperref}
\usepackage[noabbrev]{cleveref}

\newcommand{\NewTheorem}[2]{
	\newaliascnt{#1}{TheoremEnvironment}
	\newtheorem{#1}[#1]{#1}
	\aliascntresetthe{#1}
	\crefname{#1}{#1}{#2}
	\Crefname{#1}{#1}{#2}
}

\theoremstyle{definition}

\NewTheorem{Definition}{Definitions}
\NewTheorem{Remark}{Remarks}
\NewTheorem{Axiom}{Axioms}
\NewTheorem{Example}{Examples}
\NewTheorem{Observation}{Observations}
\NewTheorem{Convention}{Conventions}
\NewTheorem{Notation}{Notations}
\NewTheorem{Setting}{Settings}
\NewTheorem{Question}{Questions}
\NewTheorem{Answer}{Answers}
\NewTheorem{Conjecture}{Conjectures}
\NewTheorem{Problem}{Problems}
\NewTheorem{Solution}{Solutions}
\NewTheorem{Goal}{Goals}
\NewTheorem{Comment}{Comments}
\NewTheorem{Aim}{Aims}
\NewTheorem{Caution}{Cautions}
\NewTheorem{Exercise}{Exercises}
\NewTheorem{Examples}{Examples}
\theoremstyle{plain}
\NewTheorem{Proposition}{Propositions}
\NewTheorem{Lemma}{Lemmas}
\NewTheorem{Theorem}{Theorems}
\NewTheorem{Corollary}{Corollaries}

\crefname{enumi}{}{}
\Crefname{enumi}{}{}
\creflabelformat{enumi}{(#2#1#3)}
\crefname{enumii}{}{}
\Crefname{enumii}{}{}
\creflabelformat{enumii}{(#2#1#3)}
\crefname{enumiii}{}{}
\Crefname{enumiii}{}{}
\creflabelformat{enumiii}{(#2#1#3)}

\makeatletter
\renewcommand{\p@enumii}{}
\renewcommand{\p@enumiii}{}
\makeatother

\numberwithin{equation}{section}
\crefname{equation}{}{}
\Crefname{equation}{}{}
\creflabelformat{equation}{(#2#1#3)}

\newcommand{\SwapSymbols}[1]{
	\expandafter\let\expandafter\temporarysymbol\csname #1\endcsname
	\expandafter\let\csname #1\expandafter\endcsname\csname var#1\endcsname
	\expandafter\let\csname var#1\endcsname\temporarysymbol
}

\SwapSymbols{epsilon}
\SwapSymbols{phi}
\SwapSymbols{Gamma}
\SwapSymbols{Delta}
\SwapSymbols{Theta}
\SwapSymbols{Lambda}
\SwapSymbols{Xi}
\SwapSymbols{Pi}
\SwapSymbols{Sigma}
\SwapSymbols{Upsilon}
\SwapSymbols{Phi}
\SwapSymbols{Psi}
\SwapSymbols{Omega}

\newcommand{\cA}{\mathcal{A}}
\newcommand{\cB}{\mathcal{B}}

\newcommand{\cD}{\mathcal{D}}

\newcommand{\cF}{\mathcal{F}}

\newcommand{\cH}{\mathcal{H}}
\newcommand{\cI}{\mathcal{I}}

\newcommand{\cL}{\mathcal{L}}
\newcommand{\cM}{\mathcal{M}}

\newcommand{\cP}{\mathcal{P}}

\newcommand{\cS}{\mathcal{S}}
\newcommand{\cT}{\mathcal{T}}

\newcommand{\cW}{\mathcal{W}}

\newcommand{\bu}{^{\bullet}}

\newcommand{\To}{\longrightarrow}

\DeclareMathOperator{\Hom}{Hom}

\DeclareMathOperator{\RHom}{{\bf R}Hom}
\DeclareMathOperator{\RGamma}{{\bf R}\Gamma}
\DeclareMathOperator{\Ext}{Ext}

\DeclareMathOperator{\id}{id}

\DeclareMathOperator{\Mod}{Mod}
\let\mod\relax
\DeclareMathOperator{\mod}{mod}

\DeclareMathOperator{\Ker}{Ker}
\DeclareMathOperator{\Coker}{Coker}
\let\Im\relax
\DeclareMathOperator{\Im}{Im}

\DeclareMathOperator{\Supp}{Supp}

\title{On the category of cofinite complexes and modules}

\author{Reza Sazeedeh}
\address{Department of Mathematics, Urmia University, P.O.Box: 165, Urmia, Iran}
\email{rsazeedeh@ipm.ir and r.sazeedeh@urmia.ac.ir}

\begin{document}

\begin{abstract}
Let $A$ be a commutative noetherian ring, let $\mathfrak a$ be an ideal of $A$. In this paper, we extend Hartshorne's characterization of cofinite complexes to more general classes of rings. We also determine conditions under which Hartshorne's fourth question [H1] admits an affirmative answer. Finally, we investigate the  cofiniteness of complexes of $\frak a$-cofinite modules for rings of lower dimensions.  
\end{abstract}

\maketitle
\tableofcontents

\section{Introduction}
Throughout this paper $A$ is a commutative noetherian ring,
$\frak a$ is an ideal of $A$ and ${\bf x}=x_1,\dots x_t$ is a sequence of elements in $A$. The category of $A$-modules is denoted by $\Mod A$ and the category of finitely generated $A$-modules is denoted by $\mod A$.  An $A$-module $M$ is said to be $\frak a$-{\it cofinite} if it satisfies the following conditions:

(i) $\Supp_A M\subseteq V(\frak a);$ 

(ii) $\Ext_A^i(A/\frak a, M)$ are finitely generated for all integers $i\geq 0$.  

Hartshorne [H1] introduced the notion of an $\frak a$-cofinite module generalizing the local duality theorem for an arbitrary ideal $\frak a$ in order to answer some questions raised by Grothendieck [G]. While working towards better formulating this generalization, he raised four questions regarding $\frak a$-cofinite modules and provided answers to them in the case where $A$ is a regular ring. Recently Kawasaki [K] addressed two of these questions in the case where $A$ is a homomorphic image of a Gorenstein ring of finite Krull dimension. The primary objective of this paper is to answer to Hartshorne's questions in the context of more general rings. Our investigating also leads to an extension of the main  theorem established by Takahashi and Wakasugi [TW].  

In this paper, we deal with three full subcategories of $\Mod A$:

$\bullet$ $\cM(A,\frak a)_{\rm cof}$, the subcategory of $\Mod A$ consisting of $\frak a$-cofinite modules.

$\bullet$ $\cM^0(A,\frak a)_{\rm cof}$, the subcategory of $\Mod A$ consisting of $A$-modules satisfying the condition (ii).

$\bullet$  $\cL(A,\frak a)_{\rm cof}$, the subcategory of $\Mod A$ consisting of $A$-modules $M$ such that $H_{\frak a}^i(M)$ is $\frak a$-cofinite for all integers $i$. 

  In Section 2, we establish a criterion for cofiniteness of modules in terms of their Koszul complexes. Specifically, if ${\bf x}$ is a sequence of elements in $\frak a$ and $M$ is an $A$-module, then $M$ lies in each of the aforementioned subcategories provided its Koszul cohomology modules belong to the same subcategory.
   
In Section 3, we deal with three full subcategories of the derived category $D(A)$: 

$\bullet$ $D(A,\frak a)_{\rm cof}$, the full subcategory of $D(A)$ consisting of $\frak a$-cofinite complexes.

$\bullet$ $D_{\rm cof}(A,\frak a)$ the full subcategory of $D(A)$ consists of $A$-complexs $X\bu$ such that $H^i(X)\in \cM(A,\frak a)_{\rm cof}$ for all integers $i$.

$\bullet$ $D^0_{\rm cof}(A,\frak a)$, the full subcategory of $D(A)$ consisting of  $A$-complexs $X\bu$ such that $H^i(X\bu)\in\cM^0(A,\frak a)_{\rm cof}$ for all integers $i$.

We show that $D_{\rm cof}(A,\frak a)$ (resp. $D^0_{\rm cof}(A,\frak a)$) is a thick subcategory of $D(A)$ if and only if $\cM(A,\frak a)_{\rm cof}$ (resp. $\cM^0(A,\frak a)_{\rm cof}$)  is an abelian subcategory of $\Mod A$. If $A$ is complete with respect to $\frak a$-adic topology, we show that $D(A,\frak a)_{\rm cof}$ is a thick subcategory of $D(A)$. The following theorem originally established by Hartshorne [H1] characterizes cofinite complexes for more general rings.

\begin{Theorem}\label{hart}
Let $A$ be a ring with a dualizing complex $\cD$ and let $X\bu\in D^{+}(A)$. Consider the following conditions:\\
${\rm (1)}$ The complex $X\bu$ is $\frak a$-cofinite.\\
${\rm (2)}$ The complex $X\bu$ satisfies the following conditions:

${\rm (i)}$ $\Supp  H^i(X\bu)\subseteq V(\frak a)$ for each $i$.

${\rm (ii)}$ $\Ext^i_A(A/\frak a,X\bu)$ is finitely generated for each $i$.\\
Then ${\rm(1)}\Longrightarrow {\rm (2)}$ holds. Moreover, if $A$ is complete with respect to $\frak a$-adic topology, then ${\rm(2)}\Longrightarrow {\rm (1)}$ holds as well.
\end{Theorem}
We show that if $\cM(A,\frak a)_{\rm cof}$ is an abelian subcategory of $\Mod A$, then $\cL(A,\frak a)_{\rm cof}$ is a thick subcategory of $\Mod A$.
We prove that if $(A,\frak m)$ is a local ring and  $X\bu\in D(A,\frak a)_{\rm cof}$ is an $A$-complex in $D^+(A)$, then $\dim A/\frak a\leq\id X\bu-\inf X\bu$, where $\id X\bu$ is the injective dimension of $X\bu$. We provide characterizations for $\frak a$-cofinite  modules and complexes in terms of their Koszul complexes.

 A full abelian subcategory $\cW$ of $\Mod A$ is said to be wide if it is closed under extensions. Also a full subcategory $\cS$ of $A$-modules is said to be Serre whenever $\cS$ is closed under submodules, quotients and extensions. Let $\cH_{\frak a}(A)=\{H_{\frak a}^i(A)|\hspace{0.1cm} i\geq 0\}$ and let Sing$(A)$ be the {\it singular locus} of $A$, that is, the set of prime ideals $\frak p$ of $A$ such that the local ring $A_{\frak p}$ is not regular. One of the main results in Section 4 is  the following theorem that provides an answer to  Hartshorne's fourth question [H1] for a broader class of rings.

\begin{Theorem}\label{hartsh2}
Let $A$ admit a dualizing complex $\cD$ and let $X\bu\in D^+(A)$ be a complex of $A$-modules. If $A$ is complete with respect to $\frak a$-adic topology and $X\bu\in D_{\rm cof}(A,\frak a)$, then $X\bu\in D(A,\frak a)_{\rm cof}$. Conversely, if $X\bu\in D(A,\frak a)_{\rm cof}$ and one of the following conditions is satisfied:

 ${\rm (i)}$ $\cM^0(A,\frak a)_{\rm cof}$ is abelian,
 
${\rm (ii)}$ $\cM(A,\frak a)_{\rm cof}$ is abelian and $H_{\frak a}^i(H^i(\cD))\in\cM(A,\frak a)_{\rm cof}$ for all integers $i$, 

${\rm (iii)}$ $\cM(A,\frak a)_{\rm cof}$ is abelian, $\cH_{\frak a}(A)\subset\cM(A,\frak a)_{\rm cof}$ and {\rm Sing}$(A)\subseteq V(\frak a)$,

${\rm (iv)}$ $\cM(A,\frak a)_{\rm cof}$ is Serre and  $\cH_{\frak a}(A)\subset\cM(A,\frak a)_{\rm cof}$,

${\rm (v)}$ $A$ is a regular ring of finite Krull dimension and $\langle\cH_{\frak a}(A)\rangle_{\rm wide}\subset \cM(A,\frak a)_{\rm cof},$ where $\langle\cH_{\frak a}(A)\rangle_{\rm wide}$ is the smallest wide subcategory of $\Mod A$ containing $\cH_{\frak a}(A)$.

${\rm (vi)}$ $\sqrt{\cH_{\frak a}(A)}\subset \cM(A,\frak a)_{\rm cof}$, where $\sqrt{\cH_{\frak a}(A)}$ is the smallest Serre subcategory of $\Mod A$ containing $\cH_{\frak a}(A)$,\\
 then $X\bu\in D_{\rm cof}(A,\frak a)$.
 \end{Theorem}

Takahashi and Wakasugi [TW, Proposition 3.8] proved that $\mod A\subset\cL(A,\frak a)_{\rm cof}$ if $A$ satisfies the condition (i). In \cref{lm0}, we show that $\cL(A,\frak a)_{\rm cof}\subseteq \cM^0(A,\frak a)_{\rm cof}$. Moreover, if the condition (i) is satisfied (without any other conditions on $A$), then $\cL(A,\frak a)_{\rm cof}=\cM^0(A,\frak a)_{\rm cof}$. This  generalizes  their result as $\mod A\subset\cM^0(A,\frak a)_{\rm cof}$. We prove that if $D_{\rm cof}^0(A,\frak a)$ is a thick subcategory of $D(A)$, then the functor $\RGamma_{\frak a}:D^0_{\rm cof}(A,\frak a)\To D_{\rm cof}(A,\frak a)$ is essentially surjective. We also prove the following theorem:

\begin{Theorem}\label{twoo}
 Let $M$ be a finitely generated $A$-module. Then  the following implications hold.
 
${\rm (i)}$ If ${\rm pd}_A M<\infty$, then 
$H_{\frak a}^i(M)\in\langle\cH_{\frak a}(A)\rangle_{\rm wide}$ for each $i\geq 0$.

 ${\rm (ii)}$ $H_{\frak a}^i(M)\in\sqrt{\cH_{\frak a}(A)}$ for each $i\geq 0$. 
\end{Theorem}

As the last consequence in this section, under some additional conditions on $A$, we present the following theorem  which extends [TW, Theorem 1.2]:  
  \begin{Theorem}
    Let $A$ admit a dualizing complex $\cD$ and let $A$ be complete with respect to $\frak a$-adic topology. If $A$ satisfies one of the conditions (ii), (iii), (iv), (v) and (vi) in \cref{hartsh2}, then  $$\cL(A,\frak a)_{\rm cof}=\cM^0(A,\frak a)_{\rm cof}.$$ 
 \end{Theorem}
  
 In Section 5, we investigate $\frak a$-cofiniteness of a complexe $X\bu$ of $\frak a$-cofinite modules over rings with lower dimensions.
We show that if $\dim A\leq 3$ and $\Hom_A(A/\frak a,H^i(X))$ is finitely generated for all integers $i$, then $X\in D_{\rm cof}(A,\frak a)$. Let $A$ be a local ring such that $\dim A/({\bf x})\leq 3$ and let $X\bu$ be a bounded below  complex of $A$-modules such that $\Hom_A(A/\frak a,H^i({\bf x},X\bu))$ is finitely generated for any integer $i$. We show that $K\bu({\bf x},X\bu)\in D_{\rm cof}(A,\frak a)$.

\section{A criterion for cofiniteness of modules}

In this section, we assume that $n$ is a non negative integer. We establish a criterion for cofiniteness of modules via Koszul complexes. We start with some definitions and lemmas which are needed later.

\begin{Definition} An $A$-module $M$ is said to be $\frak a$-{\it cofinite} if it satisfies the following conditions:

(i) $\Supp_A(M)\subseteq V(\frak a);$ 

(ii) $\Ext_A^i(A/\frak a, M)$ are finitely generated for all integers $i\geq 0$. 

We denote by $\cM(A,\frak a)_{\rm cof}$, the subcategory of $\Mod A$ consisting of $\frak a$-cofinite modules. We denoted by $\cM^0(A,\frak a)_{\rm cof}$, the full subcategory of $\Mod A$ consisting of $A$-modules satisfying the condition (ii). This subcategory was originally introduced by Takahashi and Wakasugi [TW] studying the cofiniteness of local cohomology modules. 
\end{Definition}

\medskip
\begin{Definition}
We denote by $\cL(A,\frak a)_{\rm cof}$ the full subcategory of $\Mod A$ consisting of $A$-modules $M$ such that $H_{\frak a}^i(M)$ is $\frak a$-cofinite for all integers $i$. The subcategory $\cL(A,\frak a)_{\rm cof}\cap \mod A$, denoted by $C_{\frak a}(A)$, was investigated by Takahashi and Wakasugi [TW].
\end{Definition}

\medskip
\begin{Lemma}\label{change}
Let $B$ be a finitely generated $A$-algebra and let $M$ be a $B$-module. Then $\Ext_A^i(A/\frak a,M)$ is finitely generated for all $i\leq n$ if and only if $\Ext_B^i(B/\frak aB,M)$ is finitely generated for all $i\leq n$. In particular $M$ is $\frak a$-cofinite if and only if it is $\frak aB$-cofinite.
\end{Lemma}
\begin{proof}
A proof similar to that of [KS, Proposition 2.15] establishes the assertion. 
\end{proof}
\medskip
A full subcategory $\cS$ of $A$-modules is said to be {\it Serre} whenever $\cS$ is closed under  submodules, quotients and extensions.
  
\begin{Lemma}\label{mel}
Let $(T^i)$ be a connected sequence of functors between abelian categories $\cA$ and $\cB$, let $f:M\to N$ be a morphism in $\cA$ and let $\cS$ be a Serre subcategory of $\cB$. If for a certain $i$, $T^i\Coker f$ and $T^{i+1}\Ker f$ belong to $\cS$, then $\Coker T^if$ and $\Ker T^{i+1}f$ belong to $\cS$.  
\end{Lemma}
\begin{proof}
See [M1, Corollary 3.2].
\end{proof}
\medskip

\begin{Lemma}\label{p}
Let $\frak b$ be an ideal of $A$ contained in $\frak a$ and let $\Ext_A^j(A/\frak a,\Ext_A^i(A/\frak b,M))$ be finitely generated for all $i,j\leq n$. Then $\Ext_A^i(A/\frak a,M)$ is finitely generated for all $i\leq n$. 
\end{Lemma}
\begin{proof}
Let $0\To M\To E^0\stackrel{d^0}\To
E^1\stackrel{d^1}\To\dots$ be an injective resolution of $M$. Splitting it, for each $i\geq 0$, we have the short exact sequence $$0\To M^i\To E^i\To M^{i+1}\To
0$$ where $M^i=\Ker d^i$. We notice that for
each $i\geq 0$, we have
$\Ext_A^{i+1}(A/\frak a,M)\cong \Ext_A^1(A/\frak a, M^i)$ and also we have $\Ext_A^{i+1}(A/\frak b,M)\cong \Ext_A^1(A/\frak b, M^i).$
Hence, for each $i\geq 0$, there is an exact sequence
$$0\To (0:_{M^i}\frak b)\To(0:_{E^i}\frak b)\stackrel{f_i}\To(0:_{M^{i+1}}\frak b)\To\Ext_A^{i+1}(A/\frak b,M)\To 0\hspace{0.81cm}(\dag_i).$$
If we set $\bar{A}=A/\frak b$ and $\bar{\frak a}=\frak a\bar{A}$, it follows from \cref{change} that 
$\Ext_{\bar{A}}^i(\bar{A}/\bar{\frak a},\Ext_A^j(A/\frak b,M))$ is finitely generated for all $i,j\leq n$
. The assumption implies that $\Hom_{\bar{A}}(\bar{A}/\bar{\frak a},\Coker f_i)=\Hom_{\bar{A}}(\bar{A}/\bar{\frak a},\Ext_A^{i+1}(A/\frak b,M))$ is finitely generated for each $i\leq n-1$. We now show that $\Ext_{\bar{A}}^1(\bar{A}/\bar{\frak a},\Ker f_i)=\Ext_{\bar{A}}^1(\bar{A}/\bar{\frak a},(0:_{M_i}\frak b))$ is finitely generated for each $i\leq n-1$. To do this, applying the functor $\Ext_{\bar{A}}^1(\bar{A}/\bar{\frak a},-)$ to the exact sequence  $(\dag_{i-1})$ and using the assumption and the fact that $(0:_{E^{i-1}}\frak b)$ is an injective $\bar{A}$-module, it suffices to show that $\Ext_{\bar{A}}^{2}(\bar{A}/\bar{\frak a},(0:_{M^{i-1}}\frak b))$ is finitely generated. Continuing this way with the exact sequences $(\dag_k)$ for $k=i-2,\dots,0$, we have to show that $\Ext_{\bar{A}}^{i+1}(\bar{A}/\bar{\frak a},(0:_M\frak b))$ is finitely generated; which follows directly from the assumption. Now applying \cref{mel}, we deduce that $\Coker(\Hom_{\bar{A}}(\bar{A}/\bar{\frak a},f_i))\cong\Ext_A^{i+1}(A/\frak a,M)$ is finitely generated. We also observe that $\Hom_A(A/\frak a,M)\cong\Hom_A(A/\frak a,\Hom_A(A/\frak b,M))$ is finitely generated by the assumption.
\end{proof}

Suppose that $K_{\bullet}({\bf x})$ is the Koszul complex induced by ${\bf x}$. For every  $R$-module $M$, the complex $\Hom_A(K_{\bullet}({\bf x}),M)$ is denoted by 
$K^{\bullet}({\bf x},M)$ and  $i$-th cohomology of $K^{\bullet}({\bf x},M)$ is denoted by $H^i(({\bf x},M))$ for every integer $i$. In the rest of this section, we assume that ${\bf x}=x_1,\dots,x_t$ is a sequence of elements in $\frak a$.
\medskip

\begin{Proposition}\label{co}
Let $\Ext_A^i(A/\frak a,H^j({\bf x},M))$ be finitely generated for all $i,j\leq n$. Then $\Ext_A^i(A/\frak a,M)$ is finitely generated for all $i\leq n$. 
\end{Proposition}
\begin{proof}
Set $B=A[X_1,\dots,X_t]$ where $X_1,\dots,X_t$ are indeterminate. The rings homomorphoism $\phi:B\To A$ given by $\phi(X_i)=x_i$ for $i=1,\dots,t$ is surjective and $H^j(x_1,\dots,x_t,M)\cong\Ext^j_B(B/(X_1\dots,X_t),M)$. By the assumption, for all $i,j\leq n$; $$\Ext_A^i(A/\frak a,H^j(x_1,\dots,x_t,M))\cong\Ext_A^i(A/\frak a,\Ext_B^j(B/(X_1,\dots,X),M))$$ are finitely generated. We observe that $\phi$ maps $\frak aB+(X_1,\dots,X_t)$ to $\frak a$; and hence by virtue of \cref{change}, $\Ext_B^i(B/(\frak aB+(X_1,\dots,X_t)),\Ext_B^j(B/(X_1,\dots,X),M))$ is finitely generated for all $i,j\leq n$. It now follows from \cref{p} that $\Ext_B^i(B/\frak aB+(X_1,\dots,X_t),M)$ is finitely generated for all $i\leq n$; and hence using again \cref{change}, the $A$-module $\Ext_A^i(A/\frak a,M)$ is finitely generated for all $i\leq n$. 
\end{proof}

\begin{Corollary}
Let $M$ be an $A$-module such that $H^i({\bf x},M)\in\cM^0(A,\frak a)_{\rm cof}$. Then $M\in\cM^0(A,\frak a)_{\rm cof}$.
\end{Corollary}

In the rest of this section, we assume that ${\bf x}$ is a system of generators of the ideal $\frak a$.

\begin{Proposition}\label{p1}
The following conditions are equivalent.

${\rm (i)}$ $\Ext_A^i(A/\frak a,M)$ is finitely generated for all $i\leq n$.

${\rm (ii)}$ $H^i({\bf x},M)$ is finitely generated for all $i\leq n$. 
\end{Proposition}
\begin{proof}
(i)$\Rightarrow$ (ii). Consider the Koszul complex $$K^{\bullet}({\bf x},M):0\To K^0\stackrel{d^0}\To K^1\stackrel{d^1}\To\dots\stackrel{d^{t-1}}\To K^t\To 0 $$
and assume that $Z^i=\Ker d^i$, $B^i=\Im d^{i-1}$, $C^i=\Coker d^i$ and $H^i=H^i({\bf x},M)$ for each $i$. We notice that $\frak aH^i=0$, and hence $H^i\subseteq (0:_{C^i}\frak a)$ for each $i$. By the assumption $\Hom_A(A/\frak a,K^0)$ is finitely generated and $H^0=Z^0\subseteq\Hom_A(A/\frak a,K^0)$ so that $H^0$ is finitely generated. Observing the exact sequence $0\To Z^0\To K^0\To B^1\To 0$ and using the assumption, $\Ext_A^i(A/\frak a,B^1)$ is finitely generated for all $i\leq n$. Thus the exact sequence $0\To B^1\To K^1\To C^1\To 0 $ implies that $\Ext_A^i(A/\frak a,C^1)$ is finitely generated for all $i\leq n-1$. Then $H^1$ is finitely generated as it is a submodule of $\Hom_A(A/\frak a,C^1)$. Continuing this way with the exact sequences 
$$0\To H^j\To C^j\To K^{j+1}\To C^{j+1}\To o\hspace{1cm}(\dag_j).$$
we deduce that $\Ext_A^i(A/\frak a,C^j)$ is finitely generated for all $i\leq n-j$ and all $j\leq n$ and hence $H^j$ is finitely generated for all $j\leq n$. (ii)$\Rightarrow$(i) follows from \cref{co}.
\end{proof}

\begin{Corollary}
Let $M$ be an $A$-module such that $\Supp M\subseteq V(\frak a)$. Then $M$ is $\frak a$-cofinite if and only if $\Ext_A^i(A/\frak a,M)$ is finitely generated for all $i\leq t$ 
\end{Corollary}
\begin{proof}
It is clear that $H^i({\bf x},M)=0$ for all $i>t$; and hence the assertion is obtained by \cref{p1}. 
\end{proof}

\section{Characterizations of Cofinite complexes}

For any complex of $A$-modules $$X^{\bullet}=\dots\To X^{n-1}\stackrel{\partial^{n-1}}\To X^n\stackrel{\partial^{n}}\To X^{n+1}\To\dots $$  set inf$X\bu=\inf\{n\in\mathbb{Z}|\hspace{0.1cm} H^n(X\bu)\neq 0\}$, sup$X\bu=\sup\{n\in\mathbb{Z}|\hspace{0.1cm} H^n(X\bu)\neq 0\}$ and amp$X\bu=\sup X\bu-\inf X\bu$. The derived category of $A$ is denoted by $D(A)$. We use the superscript "- , +, b" to denote the subcategory of $D(A)$ consisting of homologically below, above and two-sided bounded $A$-complexes. Also,   $D_f(A)$ is denoted the full subcategory of $D(A)$ consisting of $A$-complexes with finitely generated homology modules. The symbol $\simeq$ denotes the isomorphism in $D(A).$

 \medskip
\begin{Definition}Let $\cD$ be a dualizing complex of $A$ (for the definition and the basic propertis of dualizing complex, we refer  the readers to the textbook [H2]). An $A$-complex $X\bu\in D(A)$ is said to be $\frak a$-{\it cofinite} if there exists $Y\bu \in D_f(A)$ such that $X\bu\simeq D_{\frak a}(M)$ where $D_{\frak a}(-)=\Hom_A(-,\Gamma_{\frak a}(\cD))$. We denote by $D(A,\frak a)_{\rm cof}$, the full subcategory of $D(A)$ consisting of $\frak a$-cofinite complexes. We also denote by $D_{\rm cof}(A,\frak a)$ the full subcategory of $D(A)$ consists of $A$-complexs $X\bu$ such that $H^i(X\bu)$ is $\frak a$-cofinite for all integers $i$. Furthermore, we denote by $D^0_{\rm cof}(A,\frak a)$, the full subcategory of $D(A)$ consisting of  $A$-complexs $X\bu$ such that $H^i(X\bu)\in\cM^0(A,\frak a)_{\rm cof}$ for all integers $i$.
 \end{Definition}

\medskip

A subcategory $\cT$ of $D(A)$ is said to be {\it thick} if it is closed under direct summands and for any exact triangle $X\To Y\To Z\To \Sigma X$ in $D(A)$, if two of the complexes $X,Y,Z$ are in $\cT$, then so is the third. The following proposition provides a characterizatin of the thickness of $D_{\rm cof}(A,\frak a)$.

\begin{Proposition}\label{tabl}
$D_{\rm cof}(A,\frak a)$ is a thick subcategory of $D(A)$ if and only if $\cM(A,\frak a)_{\rm cof}$ is an abelian subcategory of $\Mod A$. 
\end{Proposition}
\begin{proof}
Let $D_{\rm cof}(A,\frak a)$ be a thick subcategory of $D(A)$ and let $f:M\To N$ be a homomorphism of $A$-modules in $\cM(A,\frak a)_{\rm cof}$.  Then $f$ fits into an exact triangle in $D(A)$
$$ M\stackrel{f}\To N\To {\rm con}(f)\To \Sigma M.$$ Since $D_{\rm cof}(A,\frak a)$ is thick, we deduce that con$(f)\in D_{\rm cof}(A,\frak a)$. Hence $H^{-1}({\rm con}(f))=\Ker f$ and $H^0({\rm con}(f))=\Coker f$ are $\frak a$-cofinite modules. Conversely, assume that $\cM(A,\frak a)_{\rm cof}$ is abelian and $X\bu\To Y\bu\To Z\bu\To \Sigma X\bu $ is an exact triangle in $D(A)$ such that $X\bu,Y\bu\in D_{\rm cof}(A,\frak a)$. The long exact sequence  of $A$-modules $$\dots\To H^i(X\bu)\To H^i(Y\bu)\To H^i(Z\bu)\To H^{i+1}(X\bu)\to \dots$$ and the fact that $\cM(A,\frak a)_{\rm cof}$ is closed under extensions, imply that $Z\bu\in D_{\rm cof}(A,\frak a)$. Now, assume that  $X\bu\in D_{\rm cof}(A,\frak a)$ and $Y\bu\in D(A)$ is a direct summand of $X\bu$. For each integer $i$, the functor $H^i$ is additive so that $H^i(Y\bu)$ is a direct summand of $H^i(X\bu)$. Consequently, $Y\bu\in D_{\rm cof}(A,\frak a)$.
\end{proof}

The following proposition provides a characterizatin of the thickness of $D_{\rm cof}^0(A,\frak a)$. 
\medskip
\begin{Proposition}\label{abl2}
$D^0_{\rm cof}(A,\frak a)$ is a thick subcategory of $D(A)$ if and only if $\cM^0(A,\frak a)_{\rm cof}$ is an abelian subcategory of $\Mod A$. 
\end{Proposition}
\begin{proof}
A proof similar to that of \cref{tabl} establishes  the assertion.
\end{proof}

A full subcategory $\cT$ of $A$-Mod is said to be {\it thick} if it is closed under dirct summands, extensions, kernel of epimorphisms and cokernel of monomorphisms. 

\medskip
\begin{Proposition}
 Let $\cM(A,\frak a)_{\rm cof}$ be an abelian subcategory of $\Mod A$. Then $\cL(A,\frak a)_{\rm cof}$ is a thick subcategory of $\Mod A$.
\end{Proposition}
\begin{proof}
It is clear that $\cL(A,\frak a)_{\rm cof}$ is closed under direct summands. Let $0\To L\To M\To N\to 0$ be an exact sequence in $\Mod A$ such that $M,N\in\cL(A,\frak a)_{\rm cof}$. Thus
$\RGamma_{\frak a}(M),\RGamma_{\frak a}(N)\in D_{\rm cof}(A,\frak a)$. On the other hand, we have an exact triangle  
$L\To M\To N\To \Sigma L$ which
 gives rise to the following exact triangle in $D(A)$
 $$\RGamma_{\frak a}(L)\To\RGamma_{\frak a}(M)\To \RGamma_{\frak a}(N)\To \Sigma \RGamma_{\frak a}(L).$$
  We observe that  $D_{\rm cof}(A,\frak a)$ is a thick subcategory of $D(A)$ by \cref{tabl} and  hence $\RGamma_{\frak a}(L)\in D_{\rm cof}(A,\frak a)$. Consequently, $L\in\cL(A,\frak a)_{\rm cof}$. 
 \end{proof}

The Koszul cohomology of a complex $X^{\bullet}$ is defined as  $K^{\bullet}({\bf x},X^{\bullet})=\Hom_A(K_{\bullet}({\bf x}),X^{\bullet})$ and the $i$-th cohomology of  $K^{\bullet}({\bf x},X^{\bullet})$ is denoted by $H^i({\bf x},X^{\bullet}).$  

 \medskip
\begin{Proposition}\label{ccc}
Let $n$ be an integer, let ${\bf x}$ be a system of generators of $\frak a$, and let $X^{\bullet}$ be a complex in $D^+(A)$ such that $\Ext_A^i(A/\frak a,H^j(X))$ is finitely generated for all $i,j\leq n$. Then 

${\rm (i)}$ $H^i({\bf x},X^{\bullet})$ is finitely generated for all $i\leq n$.

 ${\rm (ii)}$ $\Ext^i(A/\frak a,X^{\bullet})$ is finitely generated for all $i\leq n$.
\end{Proposition}
\begin{proof}
Without loss of generality, we may assume that inf$X^{\bullet}=0$. (i) We first assume that $X^{\bullet}\in D^{\rm b}(A)$. We proceed by induction on amp$X^{\bullet}=m$. For the case $m=0$, assume that $\cI:0\To I^0\To I^1\To\dots$ is an injective resolution of $X$. Then $\cI$ is an injective resolution of $H^0(X\bu)$. The isomorphisms 
$K\bu({\bf x},X\bu)\simeq K\bu({\bf x},\cI)\simeq K\bu({\bf x},H^0(X\bu))$ in $D(A)$ and  \cref{p1} establish the assertion. Now, assume that $m>0$. By [H2, Lemma 7.2], there exists an exact sequence of complexes $$0\To H^0(X^{\bullet})\To\sigma'_{\geq 0}X^{\bullet}\To \sigma_{\geq 1}X^{\bullet}\To 0\hspace{0.732cm}(\dag)$$ where 
 $\sigma'_{\geq 0}(X^{\bullet})=0\To X^0/\Im \partial^{-1}\To X^{1}\stackrel{\partial^1}\To\dots$ and $\sigma_{\geq 1}(X^{\bullet})=0\To \Im \partial^0\To X^{1}\To \dots$. We observe that $X\bu\cong\sigma'_{\geq 0}(X^{\bullet})$ in $D^{\rm b}(A)$. In view of [CFH, Example 2.3.18], application of the functor $\Hom_A(K({\bf x}),-)$ to the exact sequence $(\dag)$ gives rise to the exact sequence of complexes $$0\To K^{\bullet}({\bf x},H^0(X^{\bullet}))\To  K^{\bullet}({\bf x},\sigma'_{\geq 0}(X^{\bullet}))\To  K^{\bullet}({\bf x},\sigma_{\geq 1}(X^{\bullet}))\To 0.$$ Since amp$(\sigma_{\geq 1}(X^{\bullet}))<{\rm amp}X^{\bullet}$, the induction hypothesis yields the assertion. Now assume that $X\bu\in D^+(A)$ and $K_{\bullet}({\bf x}):=0\To K^{-t}\To K^{-t+1}\To\dots \To K^ 0\To\ 0$ is the Koszul homology induced by $\bf x$. There is an exact sequence of complexes 
 $$0\To \tau_{>t+n}(X)\To X\To \tau_{\leq t+n}(X)\To 0$$ where $\tau_{\leq t+n}(X)=\dots \To X^{t+n-1}\To X^{t+n}\To 0$ and $\tau_{>t+n}(X)=0\To X^{t+n+1}\To X^{t+n+2}\To\dots$. Since $t>0$, we have $\Hom(K_{\bullet}({\bf x}),\tau_{>t+n}(X))^i=\bigoplus _{v=-t}^0\Hom(K^v,\tau_{>t+n}(X)^{v+i})=0$ for all $i\leq n$ and $\tau_{\leq t+n}(X)\in D^{\rm b}(A)$. We also observe that $H^i(X)\cong H^i(\tau_{\leq t+n}(X))$ for all $i\leq n$ as $t>0$. Thus, applying $\Hom_A(K_{\bullet}({\bf x}),-)$ and using the first case, we deduce that  $H^i({\bf x},X)\cong H^i({\bf x},\tau_{\leq n+t}(X))$ is finitely generated for all $i\leq n$. (ii) The case $X\bu\in D^{\rm b}(A)$ follows from the following exact triangle in $D(A)$ induced by $(\dag)$ 
 $$\RHom_A(A/{\frak a},H^0(X))\To \RHom_A(A/{\frak a},X)\To \RHom_A(A/{\frak a},\sigma_{\geq 1}(X))\To$$$$\Sigma\RHom_A(A/{\frak a},H^0(X))$$ and an easy induction on amp$(X)$.  For the case $X\bu\in D^+(A)$, in part (i), if we replace $K({\bf x})$ by a projective resolution of $A/\frak a$, a similar proof yields the desired assertion.  
 \end{proof}
 
  We are now ready to prove the first main theorem of this paper.
 \medskip
 \begin{proof}[Proof of \cref{hart}]
To prove (1)$\Longrightarrow$(2), assume that $Y\bu$ is a complex in $D_f(A)$ such that $X\bu\simeq D_{\frak a}(Y\bu)$. Therefore $\Supp H^i(X\bu)\subseteq V(\frak a)$ for all $i$. Without loss of generality, we may assume that $\inf X\bu=0$, $\inf \cD=m$ and $\sup \cD=n$. Furthermore,
consider the following complex of $A$-modules
 $$\sigma_{\leq n+1}:=\dots\stackrel{\partial^{n+1}}\To Y^{n}\To \ker \partial^{n+1}\To 0.$$ For every integer $k\geq 0$, we have the following equalities $$\Hom(\sigma_{\leq n+1}(Y\bu),\Gamma_{\frak a}(\cD))^k=\prod_{i\in\mathbb{Z}}\Hom(\sigma_{\leq n+1}(Y\bu)^i,\Gamma_{\frak a}(\cD^{i+k}))=\prod_{i\leq n-k}\Hom(\sigma_{\leq n+1}(Y\bu)^i,\Gamma_{\frak a}(\cD^{i+k}))$$$$=\prod_{i\leq n-k}\Hom(Y^i,\Gamma_{\frak a}(\cD^{i+k}))=\Hom(Y\bu,\Gamma_{\frak a}(\cD))^k.$$
Therefore $X\bu\simeq \Hom(Y\bu,\Gamma_{\frak a}(\cD))\simeq\Hom(\sigma_{\leq n+1}(Y),\Gamma_{\frak a}(\cD))$. As $H^i(\sigma_{\leq n+1}(Y\bu))= H^i(Y\bu)$ for each $i\leq n+1$ and $H^i(\sigma_{\leq n+1}(Y\bu))=0$ for each $i>n+1$, we find that $\sigma_{\leq n+1}(Y)\in D_f^-(A)$. Hence replacing $Y\bu$ by $\sigma_{\leq n+1}(Y\bu)$, we may assume that $Y\bu\in D_f^-(A)$. We now have the following isomorphism in $D(A)$: $$\RHom_A(A/\frak a,X)\simeq \RHom_A(Y\bu,\Hom_A(A/\frak a,\Gamma_{\frak a}(\cD)))$$$$\simeq\RHom_A(Y\bu,\Hom(A/\frak a,\cD))\simeq\RHom_A(A/\frak a\otimes^{\bf L}_AY\bu,\cD).$$ Since $Y\bu\in D_f^-(A)$, by [CFH, Proposition 2.5.19], there exists a complex $\cP:=\dots \To P^{s-1}\To P^s\To 0$ of finitely generated projective modules such that $A/\frak a\otimes^{\bf L}Y\bu\simeq \cP$ in $D(A)$. Thus $H^i(\RHom_A(A/\frak a\otimes^{\bf L}_AY\bu, H^j(\cD))\cong H^i(\Hom(\cP,H^j(\cD)))$ is finitely generated for all integers $i,j$. We prove by induction on amp$\cD=n$ that $\Ext^i_A(A/\frak a,X\bu)$ is finitely generated for all $i$. By [H2, Lemma 7.2], there is a triangle in $D(A)$ 
$$ H^m(\cD)\To \cD\To \sigma_{\geq m+1}(\cD)\To \Sigma H^m(\cD)$$ 
where $\sigma_{\geq m+1}(\cD)=0\To\Im d^m\To D^{m+1}\stackrel{d^{m+1}}\To\dots\To D^n\To 0$. We notice that ${\rm amp}\sigma_{\geq m+1}(\cD)<{\rm amp}\cD$ and  $H^i(\sigma_{\geq m+1}(\cD))=H^i(\cD)$ for all $i\geq m+1$. Applying the functor $\RHom_A(A/\frak a\otimes^{\bf L}Y\bu,-)$ to the above exact triangle and using the induction hypothesis, the assertion is required. (2)$\Longrightarrow$(1). For every $Y\bu\in D_f(A)$, it follows from the affine duality theorem [H1, Theorem 4.1] that $Y\bu\simeq D_{\frak a}D_{\frak a}(Y\bu)$. Moreover, according to [CFH, Proposition 7.6.16], we have $\Hom(Y\bu,\cD)\in D_f(A)$. Now, a similar proof as given in [H1, Theorem 5.1] establishes the assertion.   
\end{proof}

\medskip
 An immediate consequence can be given as follows.
 
\begin{Corollary}
Let $A$ be a ring with a dualizing complex which is complete with respect to $\frak a$-adic topology, and let $X\in D^0_{\rm cof}(A,\frak a)$ be an $A$-complex in $D^+(A)$. Then $\RGamma_{\frak a}(X)$ is an $\frak a$-cofinite complex.
\end{Corollary}
\begin{proof}
It follows from \cref{ccc} that $\Ext_A^i(A/\frak a,X\bu)$ is finitely generated and so the isomorphism $\RHom_A(A/\frak a,X\bu)\simeq \RHom_A(A/\frak a,\RGamma_{\frak a}(X\bu))$ in $D(A)$ and \cref{hart} conclude that $\RGamma_{\frak a}(X\bu)$ is an $\frak a$-cofinite complex.   
\end{proof}
    
\medskip
   
\begin{Proposition}
 Let $A$ be a ring with a dualizing complex which is complete with respect to $\frak a$-adic topology. Then $D(A,\frak a)_{\rm cof}$ is a thick subcategory of $D(A)$.
\end{Proposition}
\begin{proof}
If $X\bu\in D(A,\frak a)_{\rm cof}$ and $Y\bu$ is a direct summand of $X\bu$, then $\RHom_A(A/\frak a,Y\bu)$ is a direct summand of $\RHom_A(A/\frak a,X\bu)$. We also observe that $\Supp H^i(Y\bu)\subset V(\frak a)$ for all integers $i$. Then \cref{hart} implies that $Y\in D(A,\frak a)_{\rm cof}$. Let $X\bu\To Y\bu\To Z\bu\To \Sigma X\bu$ be an exact triangle in $D(A)$ such that $X,Y\in D(A,\frak a)_{\rm cof}$. Then, it follows from \cref{hart} and the exact triangle that $\Supp H^i(Z\bu)\subset V(\frak a)$ for all integers $i$. Applying the functor $\RHom_A(A/\frak a,-)$ to the above exact triangle, we have an exact triangle in $D(A)$
$$\RHom_A(A/\frak a,X\bu)\To \RHom_A(A/\frak a,Y\bu)\To \RHom_A(A/\frak a,Z\bu)\To \Sigma \RHom_A(A/\frak a,X\bu).$$
 Computing the cohomology modules and using \cref{hart}, we conclude that $Z\in D(A,\frak a)_{\rm cof}$.
\end{proof}

 \medskip

For an $A$-complex $X\bu$, the {\it injective dimension} of $X\bu$ is defined as the smallest integer $n$ for which there exists a complex $\cI:\dots\To I^{n-1}\To I^n\To 0$ of injective modules such that $X\bu\simeq \cI$. If $(A,\frak m)$ is a local ring of dimension $d$, then a result established by Roberts [R] states that for any non exact complex of injective modules $$\cI\bu:=0\To I^0\To I^1\To\dots\To I^d\To 0$$ in $D_f(A)$, we have $I^d\neq 0$. This brings us to the following supplementary result related to $D(A,\frak a)_{\rm cof}$ .
 
\medskip

\begin{Proposition}
Let $(A,\frak m)$ be a local ring and let $X\bu\in D(A,\frak a)_{\rm cof}$ be an $A$-complex in $D^+(A)$. Then $\dim A/\frak a\leq\id X\bu-\inf X\bu$, where $\id X\bu$ is the injective dimension of $X\bu$. 
\end{Proposition}
\begin{proof}
We may assume that inf$X\bu=0$. We observe that if $\id X\bu=\infty$, there is nothing to prove. Otherwise, assume that $\cI\bu:=0\To I^0\To\dots\to I^t\To 0$ is an injective resolution of $X\bu$. Then $\Hom(A/\frak a,\cI\bu)$ is a complex of injective $A/\frak a$-modules.  Since $\inf X\bu=0$, we have $H^0(X\bu)=\Ker(I^0\To I^1)=K^0\neq 0$ and it follows from \cref{hart} that $\Hom(A/\frak a,K^0)\neq 0$ as Supp$K^0\subset V(\frak a)$. This shows that $\Hom(A/\frak a,I\bu)$ is not exact and \cref{hart} implies that $H^i(\Hom(A/\frak a,\cI\bu))$ is finitely generated for all $i$. Thus [R] implies that $\Hom(A/\frak a,I^d)\neq 0$ where $d=\dim A/\frak a$. Consequently $d\leq t$.  
\end{proof}

\medskip
The following result establishes a connection between $\frak a$-cofiniteness of modules and $\frak a$-cofiniteness of their Koszul complexes.
 
\begin{Proposition}\label{prk}
Let $A$ be a ring with a dualizing complex, let ${\bf x}\in\frak a$ and let $M$ be an $A$-module such that $\Supp M\subset V(\frak a)$. If  $K\bu({\bf x},M)$ is an $\frak a$-cofinite complex, then $M$ is an $\frak a$-cofinite $A$-module. Moreover, if $A$ is complete with respect to $\frak a$-adic topology, then the converse holds as well.  
\end{Proposition}
\begin{proof}
Since ${\bf x}\in\frak a$, every differential in  the complex $A/\frak a\otimes_AK_{\bullet}({\bf x})$ is zero. Then there are the following isomorphisms in $D(A)$ $$\RHom_A(A/\frak a,K\bu({\bf x},M))\simeq \RHom_A(A/\frak a\otimes_AK_{\bullet}({\bf x}),M)$$$$=\RHom_A(\bigoplus_{{i=0}}^t\Sigma^{-i}(A/\frak a)^{\binom{t}{i}},M)\simeq \bigoplus_{{i=0}}^t\Sigma^{i}\RHom_A(A/\frak a,M)^{\binom{t}{i}}.$$ 
If $K\bu({\bf x},M)$ is $\frak a$-cofinite, \cref{hart} implies that $\Ext^i(A/\frak a,K({\bf x},M))$ is finitely generated for all $i$. Now, the above isomorphisms force that $\Ext^i(A/\frak a,M)$ is finitely generated for all $i$ so that $M$ is $\frak a$-cofinite. The converse of the implication also follows by the above isomorphisms and \cref{hart}.  
\end{proof}
A similar result can be established for $A$-complexes.
\begin{Proposition}\label{kos}
Let $A$ be a ring with a dualizing complex $\cD$, let ${\bf x}\in\frak a$ and let $X^{\bullet}\in D^{+}(A)$. If $X\bu$ is an $\frak a$-cofinite complex, then so is $K({\bf x},X\bu)$. Moreover, If $A$ is complete with respect to $\frak a$-adic topology and $\Supp H^i(X\bu)\subset V(\frak a)$ for all integers $i]$, then the converse holds.
\end{Proposition}
\begin{proof}
Assume that $X\bu$ is $\frak a$-cofinite. Then there exists $Y\bu\in D_{\rm f}(A)$ such that $X\simeq \Hom(Y,\Gamma_{\frak a}(D))$. We observe that $K\bu({\bf x})\otimes_AY\in D_f(A)$. Then there is the following isomorphism in $D(A)$ which yields the assertion 
$$K\bu({\bf x},X\bu)\simeq\Hom(K\bu({\bf x})\otimes_AY\bu,\Gamma_{\frak a}(\cD)).$$
 Conversely, similar to the proof of \cref{prk}, we have the following isomorphism in $D(A)$
$$\RHom_A(A/\frak a,K\bu({\bf x},X\bu))\simeq \bigoplus_{{i=0}}^t\Sigma^{i}\RHom_A(A/\frak a,X\bu)^{\binom{t}{i}}$$
and so the result is obtained by \cref{hart}. 
\end{proof}

\section{On a question of Hartshorne}
In this section, we aim to address a question posed by Hartshorne [H1]. We explore several  conditions under which Hartshorne's fourth question has an affirmative response.  A full abelian subcategory $\cW$ of $\Mod A$ is said to be {\it wide} if it is closed under extensions. We first provide a lemma which is crucial in our investigation. This can be a viewed as generalization of [K, Lemma 4] investigated by Kawasaki.

\medskip

\begin{Lemma}\label{kaw}
Let $\cW$ be a wide  subcategory  of $\Mod A$ and let $E_2^{p,q}\Longrightarrow  H^n$  be a convergent spectral sequence of $A$-modules such that each $E_2^{p,q}\in \cW$. Then $H^n\in\cW$  for every integer $n$.
\end{Lemma}
\begin{proof}
Since $\cW$ is wide, an easy induction on $r$ implies that $E_r^{p,q}\in \cW$ for all integers $p,q$, where  by the sequence of $A$-modules 
  $$E_{r-1}^{p-r+1,q+r}\stackrel{d^{p-r+1,q+r}}\To E_{r-1}^{p,q}\stackrel{d^{p,q}}\To E_{r-1}^{p+r-1,q-r},$$
 we have $E_r^{p,q}=\Ker d^{p,q}/\Im d^{p-r+1,q+r}$. For each integer $n$,  the $A$-module $H^n$ has a finite filtration 
 $$0=\Phi^{n+1}H^{n}\subset\dots\subset\Phi^{1}H^{n}\subset \Phi^{0}H^{n}\subset H^{n}$$ 
 where $\Phi^pH^n/\Phi^{p+1}H^n\cong E_{\infty}^{p,n-p}=E_r^{p,n-p}$ for all $p$ and  $r\gg 0$. The preceding argument  implies that each quotient $\Phi^pH^n/\Phi^{p+1}H^n$ belongs to $\cW$ for all $p$. Since  $\cW$ is closed under extensions, we deduce that $H^n\in\cW$ for all integers $n$.
\end{proof}
\medskip
\begin{Lemma}\label{prp}
Let $\cM(A,\frak a)_{\rm cof}$ be an abelian category and let $X\bu$ be an $A$-complex in $D^+(A)$ such that $H_{\frak a}^i(H^j(X\bu))$ is an $\frak a$-cofinite module for any integers $i,j$. Then $\RGamma_{\frak a}(X\bu)\in D_{\rm cof}(A,\frak a)$.
\end{Lemma}
\begin{proof}
Since $X\bu\in D^+(A)$, it has an injective resolution $\cI\bu\in D^+(A)$. Then for every $p,q$, we have $H_{\frak a}^p(H^q(X\bu))\cong H_{\frak a}^p(H^q(\cI\bu))$. If we set $E_2^{p,q}:=H_{\frak a}^p(H^q(X\bu))$, by virtue of [GM, III,7.13,7.14], the Cartan-Eilenberg resolution of $\cI$ can be lied in the first quadrant. In particular, we have a convergent spectral sequence $$E_2^{p,q}:=H_{\frak a}^p(H^q(X\bu))\Longrightarrow H^{p+q}({\bf R}\Gamma_{\frak a}(X\bu)).$$ We observe that $\cM(A,\frak a)_{\rm cof}$ is a wide subcategory of $\Mod A$; and hence it follows from \cref{kaw} that $H^j({\bf R}\Gamma_{\frak a}(X\bu))$ is an $\frak a$-cofinite module for every integer $j$.
\end{proof}
\medskip

\begin{Proposition}\label{hloc}
Let $\cM(A,\frak a)_{\rm cof}$ be an abelian category, let $A$ be a ring with a dualizing complex $\cD$ such that $H_{\frak a}^i(H^j(\cD))$ is an $\frak a$-cofinite module for any integers $i,j$ and let $X\in D^+(A)$ be an $\frak a$-cofinite complex. Then $H^i(X)$ is an $\frak a$-cofinite module for every integer $i$. 
\end{Proposition}
\begin{proof}
According to \cref{prp}, we have $\Gamma_{\frak a}(\cD)\in D_{\rm cof}(A,\frak a)$ so that $H^i(\Gamma_{\frak a}(\cD))$ is  an $\frak a$-cofinite module for every integer $i$. Let $M$ be a finitely generated $A$-module. Then $M$ has a free resolution $\cF:=\dots F_1\To F_0\To 0$ such that each $F_i$ is finitely generated. Applying $\Hom_A(-,H^q(\Gamma_{\frak a}(\cD)))$ to $\cF$ and using the fact that $\cM(A,\frak a)$ is abelian , we conclude $\Ext^p(M,H^q(\Gamma_{\frak a}(\cD)))$ is an $\frak a$-cofinite module  for any integers $p,q$. On the other hand, we have a convergent spectral sequence $$E_2^{p,q}:=\Ext^p(M,H^q(\Gamma_{\frak a}(\cD)))\Longrightarrow H^{p+q}(D_{\frak a}(M)).$$
It follows from \cref{kaw} that $H^{n}(D_{\frak a}(M))$ is $\frak a$-cofinite for all $n$. Furthermore, Since $X$ is $\frak a$-cofinite, by a similar argument as mentioned in the proof of \cref{hart}, there exists $Y\bu\in D_f^-(A)$ such that $X\simeq D_{\frak a}(Y\bu)$. Moreover, we have the following convergent spectral sequence 

$$E_2^{p,q}:=H^p(D_{\frak a}(H^q(Y\bu))\Longrightarrow H^{p+q}(D_{\frak a}(Y\bu)).$$
The previous argument and \cref{kaw} imply that $H^i(D_{\frak a}(Y\bu))$ is an $\frak a$-cofinite module for every integer $i$ so that $X\bu\in D_{\rm cof}(A,\frak a)$.  
\end{proof}

\medskip
\begin{Proposition}\label{aba}
Let $M$ be an $A$-module such that $H_{\frak a}^i(M)$ is an $\frak a$-cofinite module for every $i\geq 0$. Then $\Ext_A^i(A/\frak a,M)$ is finitely generated for all $i\geq 0$.
\end{Proposition}
\begin{proof}
See [NS, Proposition 3.4].
\end{proof}
\medskip

The following proposition generalizes [TW, Proposition 3.8] as $\mod A\subset\cM^0(A,\frak a)_{\rm cof}$. The proof is the same as in [TW, Proposition 3.8], but since some aspects of the details differ, we provide a sketch of the proof.

\begin{Proposition}\label{lm0}
There is an inclusion $\cL(A,\frak a)_{\rm cof}\subseteq \cM^0(A,\frak a)_{\rm cof}$. Moreover,
if $\cM^0(A,\frak a)_{\rm cof}$ is an abelian subcategory of $\Mod A$, then $\cL(A,\frak a)_{\rm cof}=\cM^0(A,\frak a)_{\rm cof}$.  
\end{Proposition}
\begin{proof}
The inclusion follows from \cref{aba}. Now, assume that $\cM^0(A,\frak a)_{\rm cof}$ is an abelian and $M\in \cM^0(A,\frak a)_{\rm cof}$. Let $x_1,\dots,x_n$ be elements in $\frak a$. We first prove by induction on $n$ that $H_{(x_1,\dots,x_n)}^i(M)\in\cM^0(A,\frak a)_{\rm cof}$ for all integers $i$. The case $n=0$ is clear as $H_{(0)}^0(M)=M\in \cM^0(A,\frak a)_{\rm cof}$ and $H_{(0)}^i(M)=0$ for all $i\neq 0$. If $n=1$, set $x=x_1$. Then we have an exact sequence of $A$-modules $$0\To \Gamma_{(x)}(M)\To M\To M_x\To H_{(x)}^1(M)\To 0.$$ Since $x\in\frak a$, we have $\frak aA_x=A_x$ so that $\Ext_A^i(A/\frak a,M_x)\cong \Ext_{A_x}^i(A_x/\frak aA_x,M_x)=0$ for all integers $i$. This implies that $M_x\in \cM^0(A,\frak a)_{\rm cof}$ and since $\cM^0(A,\frak a)_{\rm cof}$ is abelian, the above exact sequence forces that $ \Gamma_{(x)}(M),H_{(x)}^1(M)\in\cM^0(A,\frak a)_{\rm cof}$. We notice that $H_{(x)}^i(M)=0$ for all $i>1$. For $n\geq 2$, using the Mayer-Vietoris sequence [BS, 3.2.3], as mentioned in [TW, Proposition 3.8], we deduce that $H_{(x_1,\dots,x_n)}^i(M)\in\cM^0(A,\frak a)_{\rm cof}$ for all integers $i$. If we set $\frak a=(x_1,\dots,x_n)$, the above argument concludes that  $H_{\frak a}^i(M)$ is $\frak a$-cofinite for all integers $i$; and consequently $M\in\cL(A,\frak a)_{\rm cof}$.
\end{proof}
\medskip
\medskip
\begin{Corollary}\label{prpp}
Let $D^0_{\rm cof}(A,\frak a)$ be a thick subcategory of $D(A)$ and let 
 $X\bu\in D^0_{\rm cof}(A,\frak a)$. Then $\RGamma_{\frak a}(X)\in D_{\rm cof}(A,\frak a)$. Furthermore, the functor $\RGamma_{\frak a}:D_{\rm cof}^0(A,\frak a)\To D_{\rm cof}(A,\frak a)$ is essentially surjective.
\end{Corollary}
\begin{proof}
It follows from \cref{abl2} that $\cM^0(A,\frak a)_{\rm cof}$ is abelian and so $\cL(A,\frak a)_{\rm cof}=\cM^0(A,\frak a)_{\rm cof}$ by \cref{lm0}. Therefore, $H_{\frak a}^i(H^i(X))$ is $\frak a$-cofinite for all integers $i,j$. Now, the assertion follows from \cref{prp} as $\cM(A,\frak a)_{\rm cof}$ is abelian as well. To prove the second claim, since $\cM(A,\frak a)_{\rm cof}$ is an abelian subcategory of $\Mod A$, \cref{tabl} implies that $D_{\rm cof}(A,\frak a)$ is a thick subcategory of $D(A)$. Now, if $X\bu\in D_{\rm cof}(A,\frak a)$, by virtue of [H2, Lemma 4.6], the complex $X\bu$ has an injective resolution $\cI$ such that $\Supp I^i\subseteq V(\frak a)$ for each term $I^i$ in $\cI$. Consequently $X\bu\simeq \RGamma_{\frak a}(X\bu)$ and $X\bu\in D_{\rm cof}^0(A,\frak a)$  as $D_{\rm cof}(A,\frak a)\subset D_{\rm cof}^0(A,\frak a)$. 
\end{proof}

We denote by $\sqrt{\cS}$, the smallest Serre subcategory of $\Mod A$ containing $\cS$. For a subcategory $\cW$ of $\Mod A$, the smallest wide subcategory of $\Mod A$ containing $\cW$ is denoted by $\langle\cW\rangle_{\rm wide}$.

 For every finitely generated $A$-module $M$ and every projective resolution of $M$
$$\dots\To P_1\To P_0\To M\To 0,$$  
such that each $P_j$ is finitely generated, $\Omega_j(M)={\rm CoKer}(P_{j+1}\To P_j)$ is called the {\it $j$-th Syzygy} of $M$ with respect to this projective resolution. We observe that $M$ is the $0$-th Syzygy of every its projective resolution.  For the ideal $\frak a$, the {\it arithmetic rank} of $\frak a$, denoted by ${\rm ara}(\frak a)$, is the least number of elements in $A$ required to generate an ideal which has the same radical as $\frak a$. Set $\cH_{\frak a}(A)=\{H_{\frak a}^i(A)|\hspace{0.1cm} i\geq 0\}$. Considering the above notations, we have the following theorem.
\medskip

 \begin{Theorem}\label{twoo}
 Let $M$ be a finitely generated $A$-module. Then  the following implications hold.
 
${\rm (i)}$ If ${\rm pd}_A M<\infty$, then 
$H_{\frak a}^i(M)\in\langle\cH_{\frak a}(A)\rangle_{\rm wide}$ for each $i\geq 0$.

 ${\rm (ii)}$ $H_{\frak a}^i(M)\in\sqrt{\cH_{\frak a}(A)}$ for each $i\geq 0$. 
\end{Theorem}
 \begin{proof}
(i) We prove  by induction on pd$_AM=n$ that $H_{\frak a}^i(M)\in\langle\cH_{\frak a}(A)\rangle_{\rm wide}$ for each $i\geq 0 $; and hence the claim is obtained.
If $n=0$, then $M$ is projective. Then there exists positive integers $l,m$ such that $M\oplus A^l\cong A^m$. Then $H^i(M)=\Ker (H_{\frak a}^i(A)^m\To H_{\frak a}^i(A)^l)$ and consequently
$H^i(M)\in\langle\cH_{\frak a}(A)\rangle_{\rm wide}$ for each $i$. Now, assume that $n>0$ and so there is an exact sequence $0\To N\To P\To M\To 0$ such that $P$ is a finitely generated projective $A$-module and pd$_AN=n-1$. Fixing $i\geq 0$, we have the following exact sequence of $A$-modules  

$$\dots\To  H_{\frak a}^{i}(N)\To H_{\frak a}^{i}(P)\To H_{\frak a}^{i}(M)\To H_{\frak a}^{i+1}(N)\To H_{\frak a}^{i+1}(P)\To\dots\hspace{0.5cm}(\dag).$$
Using the induction hypothesis, the modules $H_{\frak a}^{i}(N),  H_{\frak a}^{i}(P),  H_{\frak a}^{i+1}(N)$ and $H_{\frak a}^{i+1}(P)$ belong to 
$\langle\cH_{\frak a}(A)\rangle_{\rm wide}$ so that $K^i=\Coker(H_{\frak a}^{i}(N)\To H_{\frak a}^{i}(P))$ and $L^i=\Ker( H_{\frak a}^{i+1}(N)\To H_{\frak a}^{i+1}(P))$ belong to $\langle\cH_{\frak a}(A)\rangle_{\rm wide}$. Now, since $\langle\cH_{\frak a}(A)\rangle_{\rm wide}$ is closed under extensions, we deduce that $H_{\frak a}^i(M)\in\langle\cH_{\frak a}(A)\rangle_{\rm wide}$. (ii) Suppose that $$\dots\To F_1\To F_0\To M\To 0$$ is a free resolution of $M$ such that each $F_i$ is finitely generated and $\Omega_i(M)={\rm CoKer}(F_{i+1}\To F_i)$ for each $i$. Consider the following short exact sequence of $A$-modules
$$0\To \Omega_i(M)\To F_{i-1}\To \Omega_{i-1}(M)\To 0\hspace{0.5cm}(\dag_i).$$

We first show by induction on $i\geq 1$ that $H^j_{\frak a}(\Omega_i(M))\in\sqrt{\cH_{\frak a}(A)}$ for all $j\leq i-1$. 
We notice by definition that $H_{\frak a}^j(F)\in\sqrt{\cH_{\frak a}(A)}$ for every free finitely generated $A$-module $F$ and every $j\geq 0$. Thus, for the case $i=1$, applying $\Gamma_{\frak a}(-)$ to the exact sequence $(\dag_1)$, we deduce that $\Gamma_{\frak a}(\Omega_1(M))\in\sqrt{\cH_{\frak a}(A)}$. Suppose that $i>1$. For every $j\leq i-1$, application $\Gamma_{\frak a}(-)$ to exact sequence $(\dag_i)$ gives rise the following exact sequence of $A$-modules
$$H_{\frak a}^{j-1}(\Omega_{i-1}(M))\To H_{\frak a}^{j}(\Omega_{i}(M))\To H_{\frak a}^{j}(F_{i-1}).$$ 
The induction hypothesis and the above argument conclude that $H_{\frak a}^{j-1}(\Omega_{i-1}(M))$ and $H_{\frak a}^{j}(F_{i-1})$  belong to $\sqrt{\cH_{\frak a}(A)}$ so that $\Im(H_{\frak a}^{j-1}(\Omega_{i-1}(M))\To H_{\frak a}^{j}(\Omega_{i}(M)))$ and $\Im (H_{\frak a}^{j}(\Omega_{i}(M))\To H_{\frak a}^{j}(F_{i-1})$ belong to $\sqrt{\cH_{\frak a}(A)}$ for all $j\leq i-1$. Hence, $ H_{\frak a}^{j}(\Omega_{i}(M))\in \sqrt{\cH_{\frak a}(A)}$ for all $j\leq i-1$ as $\sqrt{\cH_{\frak a}(A)}$ is closed under extensions. Now, if we set $n={\rm ara} \frak a+1$, the previous argument implies that $H_{\frak a}^j(\Omega_n(M))\in\sqrt{\cH_{\frak a}(A)}$ for all $j\leq {\rm ara}(\frak a)$ and by the basic properties of local cohomology we have $H_{\frak a}^j(\Omega_n(M))=0$ for all $j>{\rm ara}(\frak a)$ (see [BS, Theorem 3.3.1]. Therefore $H_{\frak a}^j(\Omega_n(M))\in\sqrt{\cH_{\frak a}(A)}$ for all $j\geq 0$. Now by descending induction on $i$, we prove that $H_{\frak a}^j(\Omega_i(M))\in \sqrt{\cH_{\frak a}(A)}$ for all $j\geq 0$ and $i\leq n$. The case $i=n$ has been proved. For each $j\geq 0$, application of $\Gamma_{\frak a}(-)$ to the exact sequence $(\dag_{i+1})$ gives rise to the following exact sequence of $A$-modules
$$H_{\frak a}^j(F_i)\To H_{\frak a}^j(\Omega_i(M))\To H_{\frak a}^{j+1}(\Omega_{i+1}(M)).$$ 
The previous argument and the induction hypothesis imply that $H_{\frak a}^j(F_i)$ and $H_{\frak a}^{j+1}(\Omega_{i+1}(M))$ belong to $\sqrt{\cH_{\frak a}(A)}$ for all $j\geq 0$. Since $\sqrt{\cH_{\frak a}(A)}$ is closed under  submodules, quotients and extensions, we deduce that $H_{\frak a}^j(\Omega_i(M))\in \sqrt{\cH_{\frak a}(A)}$ for all $j\geq 0$.  In particular, for $\Omega_0(M)=M$, we have $H_{\frak a}^j(M)\in\sqrt{\cH_{\frak a}(A)}$ for all $j\geq 0$. 
 \end{proof}
 
 \medskip
 
 \begin{Corollary}
 Let  $A$ be a regular ring of finite Krull dimension and let $\langle\cH_{\frak a}(A)\rangle_{\rm wide}\subset \cM(A,\frak a)_{\rm cof}$. Then $H_{\frak a}^i(M)$ is $\frak a$-cofinite for all finitely generated $A$-modules $M$ and all $i\geq 0$. 
 \end{Corollary}
 \begin{proof}
 Since $A$ is regular of finite Krull dimension, every finitely generated $A$-module is of finite projective dimension; and hence the result follows from \cref{twoo}.
 \end{proof}
\medskip
\begin{Example}[Hartshorne]
Let $k$ be a field and let $A=k[x,y][[u,v]]$ be a formal power series ring over $k$. Consider the ideal $\frak a=(u,v)$ of $A$ and the finitely generated $R$-module $M=A/(xu+yv)$. Since $A$ is a regular ring of finite Krull dimension, $M$ has finite projective dimension. We observe from [H1] that $H_{\frak a}^i(A)=0$ for $i\neq 2$  so that $H_{\frak a}^i(M)=0$ for all $i\neq 1,2$. Furthermore, we have a short exact sequence of $A$-modules $0\To A\stackrel{xu+yv.}\To A\To M\To 0$ which gives rises to the following exact sequence of $A$-modules 
 $$0\To H_{\frak a}^1(M)\To H_{\frak a}^2(A)\stackrel{xu+yv.}\To H_{\frak a}^2(A)\To H_{\frak a}^2(M)\To 0.$$ Therefore, $H_{\frak a}^i(M)\in\langle\cH_{\frak a}(A)\rangle_{\rm wide}$ for each $i\geq 0$.  
\end{Example}

 We are now ready to prove another main theorem of this paper.

  \medskip
\begin{proof}[Proof of \cref{hartsh2}]
 If $A$ is complete with respect to $\frak a$-adic topology and $H^i(X\bu)\in\cM(A,\frak a)_{\rm cof}$ for all integers $i$, it follows from \cref{ccc} that $\Ext^i(A/\frak a,X\bu)$ is finitely generated for all integers $i$. We further have $\Supp H^i(X\bu)\subset V(\frak a)$ for all $i$. Hence \cref{hart} implies that $X\bu\in D(A,\frak a)_{\rm cof}$. Conversely, for the assertion under condition (i), we have $\cL(A,\frak a)_{\rm cof}=\cM^0(A,\frak a)_{\rm cof}$ and it is clear that $\cM(A,\frak a)_{\rm cof}$ is abelian. Hence $H^i(\cD)\in\cL(A,\frak a)_{\rm cof}$ as $H^i(\cD)\in\cM^0(A,\frak a)_{\rm cof}$ for all integers $i$. Consequently, $H^i(X\bu)$ is $\frak a$-cofinite for all integers $i$ by \cref{hloc}. The assertion under condition (ii) follows directly from \cref{hloc}. The conditions (iii) and (iv) imply $H_{\frak a}^i(H^j(\cD))\in\cM(A,\frak a)_{\rm cof}$ for all integers $i,j$ by [TW, Theorem 1.2]. Hence, condition (ii) is satisfied so that the assertion is obtained in this case. Under condition (v), the complex $\cD$ is an injective resolution of $A$ and so it follows from the assumption that $H_{\frak a}^i(A)=H^i(\Gamma_{\frak a}(\cD))$ is $\frak a$-cofinite for all integers $i$. Since $\langle\cH_{\frak a}(A)\rangle_{\rm wide}$ is wide, for every finitely generated $A$-module $M$, computing  $\Ext^p_A(M,H^q(\Gamma_{\frak a}(\cD)))$ through a projective resolution of $M$ and using a proof similar that given in \cref{hloc}, we find that $\Ext_A^p(M,H^q(\Gamma_{\frak a}(\cD)))\in \langle\cH_{\frak a}(A)\rangle_{\rm wide}$ for all integers $p,q$. On the other hand, we have a convergent spectral sequence $$E_2^{p,q}:=\Ext^p_A(M,H^q(\Gamma_{\frak a}(\cD)))\Longrightarrow H^{p+q}(D_{\frak a}(M)).$$
 Hence $H^n\in\langle\cH_{\frak a}(A)\rangle_{\rm wide}$ for all integers $n$ by \cref{kaw}. On the other hand, by the assumption and a proof similar to that one given in \cref{hart}, there exists $Y\bu\in D_f^-(A)$ such that $X\simeq D_{\frak a}(Y\bu)$. Furthermore, we have the following convergent spectral sequence 

$$E_2^{p,q}:=H^p(D_{\frak a}(H^q(Y\bu))\Longrightarrow H^{p+q}(D_{\frak a}(Y\bu)).$$
The previous argument implies that $E_2^{p,q}\in \langle\cH_{\frak a}(A)\rangle_{\rm wide}$  for all integers $p,q$. Thus, by \cref{kaw}, the $A$-module $H^i(D_{\frak a}(Y\bu))$ is $\frak a$-cofinite for every integer $i$ so that $X\bu\in D_{\rm cof}(A,\frak a)$. Under condition (vi), it follows from \cref{twoo} that $H_{\frak a}^p(H^q(\cD))$  is $\frak a$-cofinite for any integers $p,q$. On the other hand, by virtue of [GM, III,7.13,7.14], there exists a convergent spectral sequence $$E_2^{p,q}:=H_{\frak a}^p(H^q(\cD))\Longrightarrow H^{p+q}(\Gamma_{\frak a}(\cD)).$$ 
Since $\sqrt{\cH_{\frak a}(A)}$ is wide, $ H^{n}(\Gamma_{\frak a}(\cD))$ is $\frak a$-cofinite  for all integers $n$ by \cref{kaw}.
 Now, by a similar argument to that given in part (v), the result  
 follows.
\end{proof}

\medskip
If ${\bf x}\in\frak a$ and $A$ admits a dualizing complex which is complete with respect to $\frak a$-adic topology, then we have the following corollary. 
 \medskip
 
\begin{Corollary}
Let $\cM(A,\frak a)_{\rm cof}$ be an abelian category and let $M$ be an $A$-module such that $H^i({\bf x},M)\in\cL(A,\frak a)_{\rm cof}$ for all integers $i$. Then $\RGamma_{\frak a}(M)$ is an $\frak a$-cofinite complex. Furthermore, if one of the conditions (i), (ii), (iii) and (iv) in \cref{hartsh2} is satisfied, then $M\in\cL(A,\frak a)_{\rm cof}$.
 \end{Corollary} 
\begin{proof}
Let $\cI$ be an injective resolution of $M$. Then $\RGamma_{\frak a}(K\bu({\bf x},M))\simeq \Gamma_{\frak a}(K\bu({\bf x},\cI)).$  For every integer $n$, we have $$\Gamma_{\frak a}(K({\bf x},\cI))^n\cong \Gamma_{\frak a}(K\bu({\bf x},\cI)^n)\cong\Gamma_{\frak a}(\bigoplus_{i=-t}^0\Hom_A(K^i,\cI^{i+n})$$$$\cong\bigoplus_{i=-t}^0\Hom_A(K^i,\Gamma_{\frak a}(\cI^{i+n}))\cong K\bu({\bf x},\RGamma_{\frak a}(M))^n$$

Then there is a convergent spectral sequence $$E_2^{p,q}:=H_{\frak a}^p(H^q(K\bu({\bf x},M)))\Longrightarrow H^{p+q}({\bf x},\RGamma_{\frak a}(M)).$$
As $E_2^{p,q}$  is $\frak a$-cofinite for all $p,q$ and $\cM(A,\frak a)_{\rm cof}$ is a wide subcategory of $\Mod A$
, it follows from \cref{kaw} that $H^{n}({\bf x},\RGamma_{\frak a}(M)))$ is an $\frak a$-cofinite module for all integers $i$. Now, \cref{hartsh2} implies that $K\bu({\bf x},\RGamma_{\frak a}(M))$ is an $\frak a$-cofinite complex and so \cref{kos} implies that $\RGamma_{\frak a}(M)$ is $\frak a$-cofinite. The second claim is clear by \cref{hartsh2}.
\end{proof}

\begin{Remark}
We remark that if $\dim A\leq 2$, then $\cM(A,\frak a)_{\rm cof}$ is abelian by [M2, Theorem 7.4] and $\cL(A,\frak a)_{\rm cof}=\cM^0(A,\frak a)_{\rm cof}$ by [M2, Theorem 7.10]. Hence the condition (ii) of \cref{hartsh2} is satisfied in this case.
\end{Remark}

The following extends [TW, Theorem 1.2].

\begin{Theorem}
    Let $A$ admit a dualizing complex $\cD$ and let $A$ be complete with respect to $\frak a$-adic topology. If $A$ satisfies one of the conditions (ii), (iii), (iv), (v) and (vi) in \cref{hartsh2},  then $$\cL(A,\frak a)_{\rm cof}=\cM^0(A,\frak a)_{\rm cof}.$$ 
 \end{Theorem}
 \begin{proof}
The inclusion $\cL(A,\frak a)_{\rm cof}\subseteq \cM^0(A,\frak a)_{\rm cof}$ follows from \cref{lm0}.  To prove the equality, let $M\in\cM^0(A,\frak a)_{\rm cof}$. We have $\RHom_A (A/\frak a, \RGamma_{\frak a}(M))\simeq \RHom_A (A/\frak a,M)$ in $D(A)$. Hence, the assumption on $M$ implies that $\Ext_A^i(A/\frak a,\RGamma_{\frak a}(M))$ is finitely generated for all integers $i$. It now follows from \cref{hart} that $\RGamma_{\frak a}(M)$ is an $\frak a$-cofinite complex. Finally \cref{hartsh2} implies that $M\in\cL(A,\frak a)_{\rm cof}$. 
 \end{proof}

\section{$\frak a$-Cofinite complexes on rings with lower dimensions}

We start this section with a result about cofiniteness of cohomology of cofinite complexes. Throughout this section, $X\bu=(X\bu,\partial)$ is a complex in $D^+(A)$ such that $X^i$ is an $\frak a$-cofinite module and $n$ is an integer number.

 If $\dim A\leq 3$ and $f:M\To N$ is a homomorphism of $\frak a$-cofinite modules, [NS, Theorem 2.8(i)] proves that $\Ker f$ and $\Coker f$ are $\frak a$-cofinite modules if and only if $\Hom_A(A/\frak a,\Coker f)$ is finitely generated. We now have the following result.  
\medskip
\begin{Proposition}\label{prcom}
Let $\dim A\leq 3$. Then the following conditions are equivalent. 

${\rm (i)}$ $\Hom_A(A/\frak a,H^i(X\bu))$ is finitely generated for all $i\leq n+1$.

${\rm (ii)}$ $H^i(X\bu)$ is $\frak a$-cofinite for all $i\leq n$.
\end{Proposition}
\begin{proof}

We may assume that $\inf X\bu=0$ and assume that $Z^i=\Ker \partial^i$, $B^i=\Im \partial^{i-1}$, $C^i=\Coker \partial ^{i-1}$ and $H^i=H^i(X\bu)$ for each $i$. Then, for each $i\geq 0$, we have the following exact sequence of modules $$0\To H^i\To C^i\To X^{i+1}\To C^{i+1}\To 0 \hspace{1cm}(\dag_i).$$  In order to prove (i)$\Rightarrow$ (ii), applying $\Hom_A(A/\frak a,-)$ to the exact sequences ($\dag_i$), we deduce that $\Hom_A(A/\frak a,C^i)$ is finitely generated for every $i\leq n+1$. Therefore [NS, Theorem 2.8(i)] implies that $C^i$ is $\frak a$-cofinite for every $i\leq n+1$. Using again [NS, Theorem 2.8(i)] and the exact sequences ($\dag_i$), we deduce that $H^i$ is $\frak a$-cofinite for every $i\leq n$. (ii)$\Rightarrow$(i). Applying $\Hom_A(A/\frak a,-)$ to the exact sequences ($\dag_i$), we deduce that $\Hom_A(A/\frak a,C^i)$ is finitely generated for every $i\leq n$ and [NS, Theorem 2.8(i)] implies  that $C^i$ is $\frak a$-cofinite for all $i\leq n+1$. Finally applying $\Hom_A(A/\frak a,-)$ to $(\dag_{n+1})$, we deduce that $\Hom_A(A/\frak a,H^{n+1})$ is finitely generated.
\end{proof}

\medskip

\begin{Corollary}\label{comp}
Let $A$ be a ring of dimension $\leq 3$ with a dualizing complex which is complete with respect to $\frak a$-adic topology. If $X\bu$ is in $D^+(A)$ such that $\Hom_A(A/\frak a,H^i(X\bu))$ is finitely generated for all integers $i$, then $X\bu$ is an $\frak a$-cofinite complex.
\end{Corollary}
\begin{proof}
It follows from \cref{prcom} that $H^i(X\bu)$ is $\frak a$-cofinite for all integer $i$. Thus $\Ext_A^i(A/\frak a,X\bu)$ is finitely generated for all integers $i$ by \cref{ccc}. Therefore, $X\bu$ is $\frak a$-cofinite by \cref{hart}.
\end{proof}

\medskip

\begin{Proposition}\label{Procomp}
Let $A$ be a local ring such that $\dim A/({\bf x})\leq 3$ and let $X\bu$ be a bounded below  complex of $A$-modules. Then the following conditions are equivalent.

${\rm (i)}$ $\Hom_A(A/\frak a,H^i({\bf x},X\bu))$ is finitely generated for all $i\leq n+1$.

${\rm (ii)}$ $H^i({\bf x},X\bu)$ is $\frak a$-cofinite for all $i\leq n$.
\end{Proposition}
\begin{proof}
It follows from [CFH, Proposition 14.3.2] that $$K^{\bullet}({\bf x},X\bu)=\Hom_A(K_{\bullet}({\bf x}),X\bu)\simeq \Sigma^{-t}(K_{\bullet}({\bf x})\otimes_AX).$$ On the other hand, by virtue of [CFH, Proposition 11.4.6], we have $({\bf x})H^i(K_{\bullet}({\bf x})\otimes_AX\bu)=0$; and hence $({\bf x})H^i({\bf x},X\bu)=0$ for all integers $i$. Since $X\bu$ is a bounded below complex and $K_{\bullet}({\bf x})$ is a bounded complex, we deduce that  $K^{\bullet}({\bf x},X\bu)$ is bounded below and so we may assume that $$K^{\bullet}({\bf x},X\bu):=0\To Y^0\stackrel{d^0}\To Y^1\stackrel{d^1}\To\dots.$$ 
For convenience, set $B=A/({\bf x})$, $B^i=\Im d^i$, $C^i=\Coker d^{i-1}$ and $H^i=H^i({\bf x},X\bu)$. Then $H^i$ is a $B$-module  for all $i\geq 0$. Since for each $i$, there exist only finitely many $X^j$ occurring in $Y^i$, we deduce that each  $Y^i$ is $\frak a$-cofinite. For each $i\geq 0$, consider the following exact sequence 
$$0\To H^i\To C^i\To Y^{i+1}\To C^{i+1}\To 0\hspace{0.8cm}(\dag_i).$$  In order to prove (i) $\Rightarrow$ (ii), we prove by induction on $n$ that $H^i$ and $C^{i+1}$ are $\frak a$-cofinite for all $i\leq n$. 
 For $n=0$, applying the functor $\Hom_A(A/\frak a,-)$ to the exact sequences $(\dag_1)$, $0\To B^1\To Y^1\To C^1\To 0$ and $0\To H^0\to Y^0\To B^1\To 0$, we deduce that $\Ext_A^i(A/\frak a,H^0)$ is finitely generated for $i\leq 2$. As $H^0$ is an $B$-module, it follows from \cref{change} that $\Ext_{B}^i(B/\frak aB,H^0)$ is finitely generated for all $i\leq 2$ and hence [NS, Corollary 2.5] implies that $H^0$ is $\frak aB$-cofinite. It now follows from \cref{change} that $H^0$ is $\frak a$-cofinite. We notice that $C^0=Y^0/\Im d^{-1}=Y^0$ is an $\frak a$-cofinite module. Thus, in view of the exact sequence $(\dag_0)$, we deduce that $C^1$ is $\frak a$-cofinite. Now assume that $n>0$ and $H^i, C^{i+1}$ are $\frak a$-cofinite for all $i\leq n-1$. Applying the functor $\Hom_A(A/\frak a,-)$ to $(\dag_{n+1})$ and using the assumption, we deduce that $\Hom_A(A/\frak a,C^{n+1})$ is finitely generated. Now, applying the functor $\Hom_A(A/\frak a,-)$ to $(\dag_n)$, using the induction hypothesis and a similar proof of the case $n=0$, we conclude that $H^n$ is $\frak a$-cofinite. (ii)$\Rightarrow$ (i). By an easy induction on $n$, we can show that $C^i$ is an $\frak a$-cofinite module for any
  $i\leq n+1$. Finally applying $\Hom_A(A/\frak a,-)$ to $(\dag_{n+1})$, we deduce that $\Hom_A(A/\frak a,H^{n+1})$ is finitely generated.
\end{proof}
\medskip

\begin{Corollary}
Let $A$ be a local ring such that $\dim A/({\bf x})\leq 3$ and let $X\bu$ be a bounded below  complex of $A$-modules such that $\Hom_A(A/\frak a,H^i({\bf x},X\bu))$ is finitely generated for any integer $i$. Then $K\bu({\bf x},X\bu)\in D_{\rm cof}(A,\frak a)$.
\end{Corollary}
\begin{proof}
The result follows from \cref{Procomp}.
\end{proof}

Melkersson [M1, Corollary 7.8] proved that if $\dim A/({\bf x})\leq 2$, then $H^i({\bf x},M)$ is $\frak a$-cofinite for all $i\geq 0$. In the following corollary, we generalize this result for $\dim A/({\bf x})\leq 3$.

\begin{Corollary}\label{dthree}
Let $A$ be a local ring such that $\dim A/({\bf x})\leq 3$ and let $M$ be an $\frak a$-cofinite module. Then $K\bu({\bf x},M)\in D_{\rm cof}(A,\frak a)$ if and only if $\Hom_A(A/\frak a,H^i({\bf x},M))$ is finitely generated for any integer $i$.
\end{Corollary}
\begin{proof}
  \cref{Procomp} yields the desired result. 
\end{proof}

\end{document}